\documentclass[12pt]{article}%
\usepackage[fleqn]{amsmath}
\usepackage{amsmath,amsfonts,amssymb,graphicx,makeidx,amsthm}
\usepackage{amscd,amsfonts,amssymb,multicol}
\usepackage{amssymb, amsmath}
\usepackage{amsfonts}
\usepackage{amssymb}
\usepackage{graphicx}%
\providecommand{\U}[1]{\protect\rule{.1in}{.1in}}
\newtheorem{theorem}{Theorem}[section]

\newtheorem{example}[theorem]{Example}

\newtheorem{remark}[theorem]{Remark}

\numberwithin{equation}{section}

\begin{document}

\title{Recovery of the rod cross section shape}
\author{Vladislav V. Kravchenko$^{1}$, Sergii M. Torba$^{2}$, Alexander O.
Vatulyan$^{3}$\\{\small $^{1}$ Departamento de Matem\'{a}ticas, Cinvestav, Unidad
Quer\'{e}taro, }\\{\small Libramiento Norponiente \#2000, Fracc. Real de Juriquilla,
Quer\'{e}taro, Qro., 76230 MEXICO}\\{\small $^{2}$ Departamento de Matem\'{a}ticas, Cinvestav, Unidad
Quer\'{e}taro, }\\{\small Libramiento Norponiente \#2000, Fracc. Real de Juriquilla,
Quer\'{e}taro, Qro., 76230 MEXICO}\\$^{3}${\small Southern Federal University, Rostov-on-Don, Russia}\\{\small e-mail: vkravchenko@math.cinvestav.edu.mx,
storba@math.cinvestav.edu.mx}}
\maketitle

\begin{abstract}
A direct method for solving the inverse problem of determining the shape of
the cross section of a rod is proposed. The method is based on Neumann series
of Bessel functions representations for solutions of Sturm-Liouville
equations. The first coefficient of the representation is sufficient for the
recovery of the unknown function. A system of linear algebraic equations for
finding this coefficient is obtained. The proposed method leads to an
efficient numerical algorithm.

\end{abstract}

\section{Introduction}

The study of the coefficient inverse problems for elastic rods is one of the
subfields of the coefficient inverse problems of elasticity theory for
inhomogeneous structures. It is a relatively new branch of the inhomogeneous
elasticity theory with applications in the theory of inhomogeneous coatings,
foundation engineering, vibrodiagnostics, seismology, bio- and nanomechanics.
The problems on the steady-state vibrations of inhomogeneous rod structures,
depending on the type of the loading, may be subdivided into three types --
those considering longitudinal, bending or torsional vibrations. Such inverse
problems for bounded structures can be regarded in the framework of the
development of the theory of inverse Sturm-Liouville problems, and possess
their own specific features. The most detailed idea of the main aspects of
inverse problems of this type can be found in monographs \cite{GladwellBook},
\cite{Vatulyan 2007}. Specific inverse problems for rods are studied in
\cite{BocharovaVatulyan2009}, \cite{Jimenez et al 2004}, \cite{Morassi 2011},
\cite{Utyashev et al 2020}, \cite{WuFricke1990}. Let us note the approach related to the
construction of nonlinear integral equations of the Uryson type, the formation
of iteration processes involving solution of a number of ill-posed problems
reducing to Fredholm integral equations of the first kind with smooth kernels,
the solutions of which are computed on the basis of the Tikhonov
regularization method \cite{Vatulyan 2007}, \cite{VatulyanYurov 2020}.
Besides, we also note the method used in \cite{Utyashev et al 2020}, which is
based on series expansions and allows one to determine the shape of the rod
cross section from a special class of functions, generated by polynomials of
order not higher than some degree. Also, we would like to cite \cite{Gao2012}, \cite{GHCh2015} and \cite{BD2002} where some purely numerical methods where used to solve inverse problems for Sturm-Liouville equations in impedance form.

In this study, we discuss a new approach for one type of inverse problems for
longitudinal vibrations of a rod with a variable cross section. Consider an
elastic rod rigidly fixed at the end $x=\pi$, the vibrations in which are
excited by the force $p(t)$, applied to the end of the rod $x=0$. Let us
assume that the elastic modulus, density and cross-sectional area are
functions of the coordinate $x$. Such problem formulations arise both for rods
made of functionally graded materials, then the variables are Young's modulus
and density, and when the variable is the cross-sectional area (for example,
for conical structures) or when there is an axisymmetric cavity inside the
rod, e.g., of a cylindrical or ellipsoidal shape.

The equation of motion for an inhomogeneous rod has the form%
\[
\frac{\partial}{\partial x}\left(  E(x)F(x)\frac{\partial u}{\partial
x}\right)  =r(x)F(x)\frac{\partial^{2}u}{\partial t^{2}}.
\]
Let us further consider the problems of steady-state oscillations assuming
that $p(t)=pe^{i\omega t}$ and $u(x,t)=u(\omega,x)e^{i\omega t}$ where
$\omega$ is the oscillation frequency. Then the equation of motion takes the
form
\[
\left(  E(x)F(x)u^{\prime}(x)\right)  ^{\prime}+\omega^{2}r(x)F(x)u(x)=0,
\]
and the boundary conditions are modelled by the equalities%
\begin{equation}\label{bndcond}
u^{\prime}(\omega,0)=-\frac{p}{EF(0)},\qquad u(\omega,\pi)=0.
\end{equation}
Suppose that there is additional information about the solution in the form of
an amplitude and frequency response:
\begin{equation}
u(\omega,0)=\widetilde{f}(\omega) \label{addcond}%
\end{equation}
at some frequencies $\omega_{1}, \omega_{2}, \ldots, \omega_{m}, \ldots$
or for $\omega\in\left[  \omega_{a},\omega_{b}\right]  $. Then it is possible
to formulate the inverse problems for recovering one of the inhomogeneous
characteristics, considering two others to be known. We focus on the situation
when $E$ and $r$ are constant, and the function $F(x)$ is unknown and is to be
determined from the additional condition (\ref{addcond}).

In the present paper we develop a method for solving this coefficient inverse
problem. The method is based on the Neumann series of Bessel functions (NSBF)
representations for solutions of Sturm-Liouville equations obtained in
\cite{KNT}, \cite{KT2018Calcolo} and on the approach initiated in
\cite{Kr2019JIIP} and developed in \cite{Kr2019MMAS InverseAxis},
\cite{DKK2019MMAS}, \cite{KrBook2020}, \cite{KKK2020JPhysConf},
\cite{KST2020IP}, \cite{KT2021 IP1}, \cite{KT2021 IP2},
\cite{Kr2022Completion}, \cite{KKC2022}, \cite{KVicente2022JMS} for solving
inverse Sturm-Liouville problems. The main idea consists in reducing the
inverse problem to the problem of finding the first coefficient, denoted by
$g_{0}(x)$, of the NSBF representation. This coefficient results to be
sufficient for recovering the unknown function $F(x)$. In the first step, from
the given information of the problem, the NSBF coefficients are computed at
the end point $x=\pi$. In the second step, a system of linear algebraic
equations for the NSBF coefficients is obtained for interior points of the
interval $(0,\pi)$. This is done by computing the multiplier constants
$\beta_{k}$, which relate certain eigenfunctions normalized at the origin with
those normalized at $x=\pi$. Finally, the NSBF coefficients for all
$x\in(0,\pi)$ and, in particular, the first NSBF coefficient $g_{0}(x)$ are
obtained by solving the system of linear algebraic equations. The function
$F(x)$ is obtained from $g_{0}(x)$ immediately (see relation (\ref{F=g0}) below).

Thus, the method is direct, and the whole algorithm reduces to dealing with
the NSBF coefficients of the solutions, which are computed from the system of
linear algebraic equations. The proposed algorithm can be easily modified to problems having more general boundary conditions than those of \eqref{bndcond}, e.g., as appearing in \cite{AMS2017}, \cite{ASX2021}.

\section{Problem of determining the shape of the cross section of the rod}

We consider the longitudinal vibrations of an elastic rod which is rigidly
fixed at the end $x=\pi$. The vibrations are generated by a periodic force
$p(t)=pe^{i\omega t}$ applied at the rod end $x=0$. Young's module $E$ and the
density of the rod $r$ are assumed to be constant, and the problem consists in
determining the shape of the cross section of the rod by measuring the
vibrations at $x=0$ at some frequencies $\omega_{1}$, $\omega_{2}$,$\ldots
$,$\omega_{m}$,$\ldots$ or for $\omega\in\left[  \omega_{a},\omega_{b}\right]
$.

The corresponding mathematical model (see, e.g., \cite[p. 72]{Vatulyan 2007})
involves the Sturm-Liouville equation%
\begin{equation}
\left(  EF(x)u^{\prime}(x)\right)  ^{\prime}+\omega^{2}rF(x)u(x)=0,
\label{rodF}%
\end{equation}
where $u(x)$ describes the vibrations of the rod and $F(x)$ the area of its
cross section. The solution $u$ is regarded as a function of the independent
variable $x$ and parameter $\omega$, so we write $u=u(\omega,x)$. The boundary
conditions have the form%
\begin{equation}
u^{\prime}(\omega,0)=-\frac{p}{EF(0)},\qquad u(\omega,\pi)=0,
\label{boundary conditions}%
\end{equation}
where $p$ and $F(0)$ are given constants different from zero. The unknown
function $F(x)>0$, $x\in(0,\pi]$, needs to be recovered from the additional
information on the solution:%
\begin{equation}
u(\omega,0)=\widetilde{f}(\omega),\qquad\omega\in\Omega, \label{addinfo}%
\end{equation}
where $\widetilde{f}$ is a given function on a set $\Omega$. The function
$F(x)$ is sought in the Sobolev class $H^{2}(0,\pi)$. This inverse problem we
will call Problem A.

Problem A is akin to the problem of recovering the Sturm-Liouville problem
from its Weyl function. Indeed, as it was shown in \cite{KT2021 IP2}, if in
(\ref{boundary conditions}) the Dirichlet condition at $\pi$ is replaced by
the Neumann condition the problem reduces to the mentioned problem of
recovering from the Weyl function given at a set of points.

In \cite{Vatulyan 2007} a method for solving Problem A was proposed, based on
a reduction of the problem to a system of nonlinear integral equations, and a
result on the uniqueness for Problem A was proved in the case when
$\Omega=\left[  \omega_{a},\omega_{b}\right]  $ does not contain resonant frequencies.

\section{Reformulation of Problem A}

Let $q\in L_{2}(0,\pi)$ be real-valued. Consider the Sturm-Liouville equation
\begin{equation}
-y^{\prime\prime}+q(x)y=\rho^{2}y,\quad x\in(0,\pi),\label{SL equation}%
\end{equation}
where $\rho\in\mathbb{C}$. The solution $y$ is regarded as a function of the
independent variable $x$ and parameter $\rho$, so we write $y=y(\rho,x)$. Let
us assume that it is subject to the boundary conditions
\begin{equation}
y^{\prime}(\rho,0)-hy(\rho,0)=c,\quad y(\rho,\pi)=0,\label{bc1}%
\end{equation}
and
\begin{equation}
y(\rho,0)=f(\rho),\label{boundary f}%
\end{equation}
where $h$ and $c$ are real constants and $f(\rho)$ is a function defined on
some set $\Omega_{\rho}$. The inverse \textbf{Problem B} consists in finding
the potential $q(x)$ and the constant $h$ by the function $f(\rho)$ and the
constant $c$.

It is easy to see that Problems A and B are equivalent. The solutions of
equations (\ref{rodF}) and (\ref{SL equation}) are related as
\[
y(\rho,x)=a(x)u(\omega,x),
\]
where
\begin{equation}
a(x):=\sqrt{F(x)},\quad q(x)=\frac{a^{\prime\prime}(x)}{a(x)},\quad\rho
=\omega\sqrt{\frac{r}{E}}, \label{relation coefficients}%
\end{equation}
and the magnitudes in the boundary conditions are related by the equalities%
\begin{equation}
f(\rho)=a(0)\widetilde{f}(\omega),\quad h=\frac{a^{\prime}(0)}{a(0)},\quad
c=-\frac{p}{Ea(0)}. \label{relation boundary}%
\end{equation}

Thus, instead of solving Problem A one may solve Problem B for $f(\rho)$ and
$c$ defined by the corresponding equalities in (\ref{relation boundary}) and
subsequently calculate $a(x)=\sqrt{F(x)}$ by solving the Cauchy problem%
\begin{equation}
a^{\prime\prime}-q(x)a=0, \label{SL0}%
\end{equation}%
\begin{equation}
a(0)=\sqrt{F(0)},\quad a^{\prime}(0)=h\sqrt{F(0)}. \label{init cond a}%
\end{equation}
As we will see below, the method proposed in this paper does not require
solving this Cauchy problem, we recover directly the function $a(x)$.

\section{Solution of Problem B}

\subsection{Representations of solutions}

By $\varphi(\rho,x)$ and $S(\rho,x)$ we denote the solutions of
\eqref{SL equation} satisfying the initial conditions at the origin
\begin{equation}
\varphi(\rho,0)=1\quad\text{and}\quad\varphi^{\prime}(\rho
,0)=h\label{init cond}%
\end{equation}
and
\begin{equation}
S(\rho,0)=0\quad\text{and}\quad S^{\prime}(\rho,0)=1,\label{init cond2}%
\end{equation}
respectively, and by $T(\rho,x)$ the solution of \eqref{SL equation}
satisfying the initial conditions at $x=\pi$,
\begin{equation}
T(\rho,\pi)=0\quad\text{and}\quad T^{\prime}(\rho,\pi)=1.\label{init cond pi}%
\end{equation}

Let $\left\{  \mu_{k}^{2}\right\}  _{k=0}^{\infty}$ be the spectrum of the
Sturm-Liouville problem
\begin{align}
-y^{\prime\prime}+q(x)y &  =\rho^{2}y,\quad x\in(0,\pi),\label{SL}\\
y^{\prime}(0)-hy(0) &  =0,\quad y(\pi)=0.\label{SLboundary}%
\end{align}
Then if $\rho^{2}$ is not an eigenvalue of this problem, i.e., $\rho^{2}%
\neq\mu_{k}^{2}$, there exists a unique solution $\Phi(\rho,x)$ of the
boundary value problem (\ref{SL equation}), (\ref{bc1}). It is easy to see
that
\[
\Phi(\rho,x)=cS(\rho,x)-c\frac{S(\rho,\pi)}{\varphi(\rho,\pi)}\varphi(\rho,x).
\]
If additionally $\Phi(\rho,x)$ satisfies the condition (\ref{boundary f}),
then necessarily the equality holds
\begin{equation}
f(\rho)\varphi(\rho,\pi)+cS(\rho,\pi)=0,\text{ }\rho^{2}\neq\mu_{k}%
^{2}.\label{eqsys1}%
\end{equation}
On the other hand, we have
\begin{equation}
\varphi(\mu_{k},\pi)=0.\label{eqsys2}%
\end{equation}

We will use the following series representations for the solutions
$\varphi(\rho,x)$ and $S(\rho,x)$ obtained in \cite{KNT}.

\begin{theorem}
[\cite{KNT}]\label{Th NSBF} The solutions $\varphi(\rho,x)$ and $S(\rho,x)$
and their derivatives with respect to $x$ admit the following series
representations
\begin{align}
\varphi(\rho,x) &  =\cos\rho x+\sum_{n=0}^{\infty}(-1)^{n}g_{n}(x)j_{2n}(\rho
x),\label{phi series}\\
S(\rho,x) &  =\frac{\sin\rho x}{\rho}+\frac{1}{\rho}\sum_{n=0}^{\infty
}(-1)^{n}s_{n}(x)j_{2n+1}(\rho x),\label{S}\\
\varphi^{\prime}(\rho,x) &  =-\rho\sin\rho x+\left(  h+\frac{1}{2}\int_{0}%
^{x}q(s)\,ds\right)  \cos\rho x+\sum_{n=0}^{\infty}(-1)^{n}\widetilde{g}%
_{n}(x)j_{2n}(\rho x),\\
S^{\prime}(\rho,x) &  =\cos\rho x+\frac{1}{2\rho}\left(  \int_{0}%
^{x}q(s)\,ds\right)  \sin\rho x+\frac{1}{\rho}\sum_{n=0}^{\infty}%
(-1)^{n}\widetilde{s}_{n}(x)j_{2n+1}(\rho x),\label{Sprime}%
\end{align}
where $j_{k}(z)$ stands for the spherical Bessel function of order $k$ (see,
e.g., \cite[p. 437]{AbramowitzStegunSpF}). For every $\rho\in\mathbb{C}$ all
the series converge pointwise. For every $x\in\left[  0,\pi\right]  $ the
series converge uniformly on any compact set of the complex plane of the
variable $\rho$, and the remainders of their partial sums admit estimates
independent of $\operatorname{Re}\rho$.

The first coefficients $g_{0}(x)$, $s_{0}(x)$, $\widetilde{g}_{0}(x)$ and
$\widetilde{s}_{0}(x)$ have the form
\begin{align}
g_{0}(x) &  =\varphi(0,x)-1,\quad s_{0}(x)=3\left(  \frac{S(0,x)}{x}-1\right)
,\label{beta0}\\
\widetilde{g}_{0}(x) &  =g_{0}^{\prime}(x)-h-\frac{1}{2}\int_{0}%
^{x}q(s)\,ds,\quad\widetilde{s}_{0}(x)=\frac{s_{0}(x)}{x}+s_{0}^{\prime
}(x)-\frac{3}{2}\int_{0}^{x}q(s)\,ds,\nonumber
\end{align}
and the rest of the coefficients can be calculated following a simple
recurrent integration procedure.
\end{theorem}

\begin{remark}
Analogous series representations are valid, of course, for the solution
$T(\rho,x)$. In particular,%
\begin{equation}
T(\rho,x)=\frac{\sin\left(  \rho\left(  x-\pi\right)  \right)  }{\rho}%
+\frac{1}{\rho}\sum_{n=0}^{\infty}(-1)^{n}t_{n}(x)j_{2n+1}\left(  \rho\left(
x-\pi\right)  \right)  .\label{s=}%
\end{equation}

\end{remark}

\begin{remark}
\label{Rem q from beta0}The first equality in \eqref{beta0} shows that the
potential $q(x)$ can be recovered from $g_{0}(x)$. Indeed,
\begin{equation}
q(x)=\frac{\varphi^{\prime\prime}(0,x)}{\varphi(0,x)}=\frac{g_{0}%
^{\prime\prime}(x)}{g_{0}(x)+1}. \label{q from beta0}%
\end{equation}
The constant $h$ is also recovered directly from $g_{0}(x)$, since%
\[
h=g_{0}^{\prime}(0).
\]
Even more interesting is the relation between the function $a(x)$ from
(\ref{relation coefficients}) and $g_{0}(x)$%
\begin{equation}
a(x)=\sqrt{F(0)}\left(  g_{0}(x)+1\right)  . \label{a=g0}%
\end{equation}
Indeed, comparing the initial conditions fulfilled by $a(x)$ and
$\varphi(0,x)$ (equalities (\ref{init cond a}) and (\ref{init cond})) and
taking into account that both functions are solutions of the same equation
(\ref{SL0}) we conclude that $a(x)=\sqrt{F(0)}\varphi(0,x)$. Now the first
equality in \eqref{beta0} leads to (\ref{a=g0}).

Hence the cross section area $F(x)$ can be recovered directly from $g_{0}(x)$
by the formula%
\begin{equation}
F(x)=F(0)\left(  g_{0}(x)+1\right)  ^{2}. \label{F=g0}%
\end{equation}

\end{remark}

\begin{remark}
Notice that due to the assumption $F(x)>0$ for all $x\in\left[  0,\pi\right]
$, we have that $g_{0}(x)\neq-1$, $x\in\left[  0,\pi\right]  $. Moreover,
$\rho=0$ is not an eigenvalue of the problem (\ref{SL}), (\ref{SLboundary})
because otherwise we would have $\varphi(0,\pi)=0$ and hence $g_{0}(\pi)=-1$.
\end{remark}

\subsection{Linear algebraic systems for coefficients at the end points}

Equation (\ref{eqsys1}) gives us a simple way to compute the coefficients
$\left\{  g_{n}(\pi)\right\}  $ and $\left\{  s_{n}(\pi)\right\}  $ when
$f(\rho)$ and $c$ are known. Indeed, substitution of (\ref{phi series}) and
(\ref{S}) into (\ref{eqsys1}) leads to the equality%
\begin{equation}
f(\rho)\sum_{n=0}^{\infty}(-1)^{n}g_{n}(\pi)j_{2n}(\rho\pi)+\frac{c}{\rho}%
\sum_{n=0}^{\infty}(-1)^{n}s_{n}(\pi)j_{2n+1}(\rho\pi)=-f(\rho)\cos\rho
\pi-c\frac{\sin\rho\pi}{\rho},\label{eqsyst1}%
\end{equation}
which is valid for all $\rho^{2}\neq\mu_{k}^{2}$. Note that at $\rho=0$ from
(\ref{eqsyst1}) we obtain
\begin{equation}
f(0)g_{0}(\pi)+\frac{\pi c}{3}s_{0}(\pi)=-f(0)-\pi c.\label{eqsyst2}%
\end{equation}

When $\rho^{2}=\mu_{k}^{2}$ and hence $f(\rho)=\infty$, equation
(\ref{eqsys2}) leads to the equality
\begin{equation}
\sum_{n=0}^{\infty}(-1)^{n}g_{n}(\pi)j_{2n}(\rho\pi)=-\cos\rho\pi
.\label{eqsyst3}%
\end{equation}
Thus, the resonant frequencies, which can be relatively easily detected in
practice are useful for computing the coefficients $\left\{  g_{n}%
(\pi)\right\}  $. However, for computing the coefficients $\left\{  s_{n}%
(\pi)\right\}  $ we need to dispose of the values of $f(\rho)$ for $\rho
^{2}\neq\mu_{k}^{2}$.

Substitution of a sufficient number of values of $\rho\in\Omega_{\rho}$ and
corresponding values of $f(\rho)$ into respective equations (\ref{eqsyst1}%
)-(\ref{eqsyst3}) gives us a system of linear algebraic equations for
computing the unknown coefficients $\left\{  g_{n}(\pi)\right\}  $ and
$\left\{  s_{n}(\pi)\right\}  $. Similarly to \cite[Theorem 3.1]{KT2021 IP2}, uniqueness of the solution can be obtained in a suitable $l_2$ spaces.

Having computed the coefficients $\left\{  g_{n}(\pi)\right\}  $ and $\left\{
s_{n}(\pi)\right\}  $ allows us to compute $\varphi(\rho,\pi)$ and $S(\rho
,\pi)$, in principle, for all $\rho$. However, for solving the inverse problem
we need to transform this information on the solutions at the end point into
information on them or at least on the coefficient $g_{0}(x)$ on the whole
interval $(0,\pi)$. To this purpose the following strategy is proposed.

With the aid of the coefficients $\left\{  g_{n}(\pi)\right\}  $ and $\left\{
s_{n}(\pi)\right\}  $ compute $\left\{  \mu_{k}\right\}  $ (the square roots
of the eigenvalues of the problem (\ref{SL}), (\ref{SLboundary})). The numbers
$\mu_{k}$ are computed as zeros of the function
\begin{equation}
\varphi(\rho,\pi)=\cos\rho\pi+\sum_{n=0}^{\infty}(-1)^{n}g_{n}(\pi)j_{2n}%
(\rho\pi).\label{phi at pi}%
\end{equation}

Consider the third solution $T(\rho,x)$ satisfying the initial conditions
(\ref{init cond pi}) at $x=\pi$. Observe that for $\rho=\mu_{k}$ the solutions
$T(\mu_{k},x)$ and $\varphi(\mu_{k},x)$ become linearly dependent, so there
exist such constants $\beta_{k}\neq0$ that
\begin{equation}
T(\mu_{k},x)=\beta_{k}\varphi(\mu_{k},x).\label{relation s and phi}%
\end{equation}
Since $\varphi(\mu_{k},0)=1$ we have that
\[
\beta_{k}=T(\mu_{k},0).
\]

Note that there is a simple relation between the values of the solutions
$S(\rho,x)$ and $T(\rho,x)$ at the opposite endpoints. Namely, for all
$\rho\in\mathbb{C}$, the equality holds
\begin{equation}
S(\rho,\pi)=-T(\rho,0).\label{S=-T}%
\end{equation}
To verify it, consider the Wronskian $W(x):=S(\rho,x)T^{\prime}(\rho
,x)-S^{\prime}(\rho,x)T(\rho,x)$. It is independent of $x$ (see, e.g.,
\cite[p. 327]{Hartman}). Evaluating it at zero and taking into account
(\ref{init cond2}) we obtain $W(0)=-T(\rho,0)$, while considering it at
$x=\pi$ and using (\ref{init cond pi}) we obtain $W(\pi)=S(\rho,\pi)$. Thus we
have (\ref{S=-T}).

Hence,%
\begin{equation}
\beta_{k}=-S(\mu_{k},\pi)\label{beta_k=}%
\end{equation}
and can be computed by using the coefficients $\left\{  s_{n}(\pi)\right\}  $,
computed in the first step.

Finally, relation (\ref{relation s and phi}) gives us a system of linear
algebraic equations for computing the coefficients $\left\{  g_{n}(x)\right\}
$, $\left\{  t_{n}(x)\right\}  $ and in particular the coefficient $g_{0}(x)$.
The system has the form
\[
\sum_{n=0}^{\infty}(-1)^{n}g_{n}(x)j_{2n}(\mu_{k}x)-\frac{1}{\beta_{k}\mu_{k}%
}\sum_{n=0}^{\infty}(-1)^{n}t_{n}(x)j_{2n+1}\left(  \mu_{k}\left(
x-\pi\right)  \right)  =-\cos\mu_{k}x+\frac{\sin\left(  \mu_{k}\left(
x-\pi\right)  \right)  }{\beta_{k}\mu_{k}}%
\]
for all $x\in\left[  0,\pi\right]  $.

\section{Algorithm for solving Problem A}\label{Section5}

The numerical algorithm consists of the following steps.

\begin{enumerate}
\item Given the values of the function $\widetilde{f}(\omega)$ at a number of
points $\omega_{l}$, $l=1,\ldots,L$. Compute
\[
\rho_{l}=\omega_{l}\sqrt{\frac{r}{E}}\text{\quad and\quad}f(\rho_{l}%
)=\sqrt{F(0)}\widetilde{f}(\omega_{l}),\text{\quad}l=1,\ldots,L.
\]

\item The set $\Omega_{\rho}=\left\{  \rho_{l}\right\}  _{l=1}^{L}$ may
contain some eigenvalues of the problem (\ref{SL}), (\ref{SLboundary}), for
which we have that $f(\rho_{l})=\infty$. Reorder the set $\Omega_{\rho}$ in
such a way that $\Omega_{\rho}=\Omega_{\rho,1}\cup\Omega_{\rho,2}$ and
$\Omega_{\rho,1}=\left\{  \rho_{l}\right\}  _{l=1}^{L_{1}}$ consists of
regular values of the problem (\ref{SL}), (\ref{SLboundary}), while
$\Omega_{\rho,2}=\left\{  \rho_{l}\right\}  _{l=L_{1}+1}^{L}$ consists of the
eigenvalues of the problem (\ref{SL}), (\ref{SLboundary}).

Notice that $L_{1}$ can be equal to $L$ (if no eigenvalue of the problem
(\ref{SL}), (\ref{SLboundary}) belongs to $\Omega_{\rho}$), but $L_{1}$ must
be greater than zero, because otherwise the problem is not uniquely solvable.

The subset $\Omega_{\rho,1}$ may contain $\rho_{l}=0$.

\item \label{AlgStep3} Compute a number of the coefficients $\left\{  g_{n}(\pi)\right\}
_{n=0}^{N_{1}}$ and $\left\{  s_{n}(\pi)\right\}  _{n=0}^{N_{2}}$ from the
system of equations%
\begin{equation}
f(\rho_{l})\sum_{n=0}^{N_{1}}(-1)^{n}g_{n}(\pi)j_{2n}(\rho_{l}\pi)+\frac
{c}{\rho_{l}}\sum_{n=0}^{N_{2}}(-1)^{n}s_{n}(\pi)j_{2n+1}(\rho_{l}\pi
)=-f(\rho_{l})\cos\rho_{l}\pi-c\frac{\sin\rho_{l}\pi}{\rho_{l}}%
,\label{equations 1}%
\end{equation}
for $l=1,\ldots,L_{1}$ and if $\rho_{l}\neq0$,
\begin{equation}
f(0)g_{0}(\pi)+\frac{\pi c}{3}s_{0}(\pi)=-f(0)-\pi c,\quad\text{if }\rho
_{l}=0,\label{equations 2}%
\end{equation}%
\begin{equation}
\sum_{n=0}^{N_{1}}(-1)^{n}g_{n}(\pi)j_{2n}(\rho_{l}\pi)=-\cos\rho_{l}\pi,\quad
l=L_{1}+1,\ldots,L.\label{equations 3}%
\end{equation}
The system in this step can be overdetermined, so that the number of the
unknowns $N_{1}+N_{2}+2$ can be much smaller than the number of equations. It
is important that the number $N_{2}+1$ of the unknown coefficients $s_{n}%
(\pi)$ (which in principle can be different from the number of the unknowns
$g_{n}(\pi)$) be not greater than the number of equations of types
(\ref{equations 1}) and (\ref{equations 2}).

\item\label{AlgStep4} Using $\left\{  g_{n}(\pi)\right\}  _{n=0}^{N_{1}}$ and $\left\{
s_{n}(\pi)\right\}  _{n=0}^{N_{2}}$ compute the numbers $\left\{  \mu
_{k}\right\}  _{k=0}^{M}$ and $\left\{  \beta_{k}\right\}  _{k=0}^{M}$. The
first of them are computed as zeros of the function
\[
\varphi_{N_{1}}(\rho,\pi)=\cos\rho\pi+\sum_{n=0}^{N_{1}}(-1)^{n}g_{n}%
(\pi)j_{2n}(\rho\pi),
\]
considered on a sufficiently large interval in $\rho$. The second sequence
$\left\{  \beta_{k}\right\}  _{k=0}^{M}$, according to (\ref{beta_k=}), are
computed as
\[
\beta_{k}=-\frac{\sin\left(  \pi\mu_{k}\right)  }{\mu_{k}}-\frac{1}{\mu_{k}%
}\sum_{n=0}^{N_{2}}(-1)^{n}s_{n}(\pi)j_{2n+1}\left(  \pi\mu_{k}\right)  ,\quad
k=0,\ldots,M.
\]

\item\label{AlgStep5} Compute $g_{0}(x)$ by solving the system
\begin{align*}
&  \sum_{n=0}^{\widetilde{N}_{1}}(-1)^{n}g_{n}(x)j_{2n}(\mu_{k}x)-\frac
{1}{\beta_{k}\mu_{k}}\sum_{n=0}^{\widetilde{N}_{2}}(-1)^{n}t_{n}%
(x)j_{2n+1}\left(  \mu_{k}\left(  x-\pi\right)  \right)  \\
&  =-\cos\mu_{k}x+\frac{\sin\left(  \mu_{k}\left(  x-\pi\right)  \right)
}{\beta_{k}\mu_{k}},\quad k=0,\ldots,M.
\end{align*}
for all $x\in\left[  0,\pi\right]  $. Here $\widetilde{N}_{1}$ and
$\widetilde{N}_{2}$ can be chosen different from $N_{1}$ and $N_{2}$ above.

\item Compute $F(x)$ from (\ref{F=g0}).
\end{enumerate}

Below we discuss some numerical implementation details of the proposed algorithm. We do not provide the exact receipts, but even the simplest implementation of the proposed ideas was sufficient in our numerical tests.

Also we would like to mention that the performance of the proposed method greatly depends on the choice of the points $\omega_l$. If all the values $\omega_l$ belongs to a small interval, the corresponding numerical problem is ill-conditioned. The situation is similar to those arising for the inverse spectral problem with variable boundary parameter \cite[Section 3.14]{Chadan et al 1997}. Also, as we illustrate in Example \ref{ExExp}, increasing the density of the points $\omega_l$ has little effect. The best results are obtained when the transformed set $\rho_l$ is close to the set consisting of integers and half-integers (like a union of eigenvalues of Neumann-Neumann and Neumann-Dirichlet spectral problems).

Note that, once Step \ref{AlgStep3} is complete, the original data $\{\rho_l, f(\rho_l)\}_{l=1}^L$ is no longer used. The coefficients $\left\{  g_{n}(\pi)\right\}  _{n=0}^{N_{1}}$ and $\left\{
s_{n}(\pi)\right\}  _{n=0}^{N_{2}}$ contain all the information used for recovery of the cross section area $F$.

Step \ref{AlgStep3} requires the choice of parameters $N_1$ and $N_2$ for system \eqref{equations 1}--\eqref{equations 3}. Numerical experiments show that the condition number of the coefficient matrix grows as $N_1+N_2$ increases. However, the absolute values of the coefficients $g_n(\pi)$ and $s_n(\pi)$ rapidly decrease as $n\to\infty$, see \cite[Proposition A.3]{KT2018Calcolo}, hence one does not need many of them. Moreover, there was no significant difference in the observed decrease rate of the coefficients $g_n(\pi)$ and $s_n(\pi)$, so we can take $N_1= N_2$ for simplicity. Also, assume that no one $\rho_l=0$. We look for least squares solution of the overdetermined system \eqref{equations 1}--\eqref{equations 3}, so let us consider a discrepancy function
\begin{equation}\label{Eq Discrepancy}
\begin{split}
    Q(N, g_0,\ldots,g_N,s_0,\ldots,s_N) &= \sum_{l=1}^{L_1} \left| f(\rho_l) \sum_{n=0}^N (-1)^n g_n j_{2n}(\rho_l \pi)+f(\rho_l) \cos \rho_l\pi \right.\\
    & \left.+ \frac c{\rho_l}\sum_{n=0}^N (-1)^n s_n j_{2n+1}(\rho_l \pi)  + c\frac{\sin \rho_l \pi}{\rho_l}\right|^2\\
    &+ \sum_{l=L_1+1}^L \left| \sum_{n=0}^N (-1)^n g_n j_{2n}(\rho_l \pi) + \cos \rho_l\pi\right|^2.
\end{split}
\end{equation}

Let
\begin{equation}\label{FunctionalOrig}
Q_N = \left(\min_{ \{g_n\}_{n=0}^N, \{ s_n \}_{n=0}^N } Q(N, g_0,\ldots g_N, s_0, \ldots, s_N)\right)^{1/2}
\end{equation}
and
\begin{equation}\label{Eq argmin}
    \{ g_0^N, \ldots, g_N^N, s_0^N, \ldots, s_N^N\} = \operatorname{arg\,min}_{\{g_n\}_{n=0}^N, \{ s_n \}_{n=0}^N} Q(N, g_0,\ldots, g_N, s_0, \ldots, s_N)
\end{equation}
(for numerical solution of \eqref{Eq argmin}, we utilize the Moore-Penrose pseudoinverse).

The exact value of $Q_N$ does not increase as $N\to\infty$. Due to numerical errors and limited machine precision, its value almost stabilizes at some $N=N_0$. The value $N_0$ and the corresponding set of coefficients $\{g_0^{N_0},\ldots, g_{N_0}^{N_0}, s_0^{N_0},\ldots, s_{N_0}^{N_0}\}$ can be taken as the result of Step \ref{AlgStep3}. However, for some numerical examples (see Example \ref{Ex2Holes}) resulted coefficients $\{g_0^{N_0},\ldots, g_{N_0}^{N_0}, s_0^{N_0},\ldots, s_{N_0}^{N_0}\}$ are rather large and change significantly if one considers $N=N_0+1$ instead. Taking larger value of $N$ leads to smaller by absolute value coefficients $\{g_n^N, s_n^N\}_{n=0}^N$, with smaller variation with respect to $N$ and better overall solution of the inverse problem. For that reason we consider the functional
\begin{equation}\label{Functional opt}
    R_N=Q_N + \alpha \left( \sum_{n=0}^{N-1} \left( |g_n^N - g_n^{N-1}|^2 + |s_n^N - s_n^{N-1}|^2\right) + |g_N^N|^2 + |s_N^N|^2\right)^{1/2},
\end{equation}
where $\alpha$ is a small parameter (we take $10^{-3}$), and use this functional to obtain an optimal value for $N_1$ and $N_2$ in Step \ref{AlgStep3}.

Note that the process described in the Steps \ref{AlgStep3} and \ref{AlgStep4} is similar to the process of spectrum completion proposed in \cite{Kr2022Completion}, and we refer the reader to this article to see the numerical performance.

Since we can compute arbitrary large set of spectral data on Step \ref{AlgStep4}, there is no need to consider a functional similar to \eqref{Functional opt} on Step \ref{AlgStep5}. Good results were obtained taking $\tilde N_1$ and $\tilde N_2$ separately and equal to numbers of singular values being greater than $10^{-2}$ of the corresponding coefficient matrices. Matlab function \texttt{rank} with a parameter \texttt{tol} readily provides these numbers.

\section{Numerical examples}
The proposed algorithm can be implemented directly taking into account the comments at the end of Section \ref{Section5}. All the computations were performed in Matlab 2017 in machine precision arithmetics. For all the examples, on Step \ref{AlgStep4} we computed 1000 eigendata pairs.

\begin{example}\label{ExParabolic}
Consider a rod having quartic cross section area, i.e.,
\[
F(x) = (a+bx)^4,
\]
where $a$ and $b$ are some constants. For this cross section profile we have
\[
a(x) = (a+bx)^2,\qquad q(x) = \frac{2b^2}{(a+bx)^2},\qquad h = \frac{2b}a,
\]
and the corresponding series \eqref{phi series} and \eqref{S} have at most two non-zero coefficients $g_n$ and $s_n$, which can be seen from  \cite[Subsection 3.3]{KrT2012} and \cite[(3.3)]{KNT}. Indeed,
\begin{align}
    \varphi(\rho,x)&=\cos \rho x+\frac{bx(2a+bx)}{a^2}j_0(\rho x) +\frac{b^3x^3}{a^2(a+bx)}j_2(\rho x),\label{Ex1phi}\\
    S(\rho,x)&= \frac{\sin \rho x}{\rho} + \frac{1}{\rho}\frac{b^2x^2}{a(a+bx)} j_1(\rho x).\label{Ex1S}
\end{align}
Hence, recalling equation \eqref{eqsys1}, we obtain the explicit expression for the direct data
\begin{equation}\label{f via phi and S}
f(\rho) = -c \frac{S(\rho,\pi)}{\varphi(\rho,\pi)}.
\end{equation}

\begin{figure}[htb!]
\centering
\includegraphics[bb=0 0 432 158, width=6in,height=2.2in]
{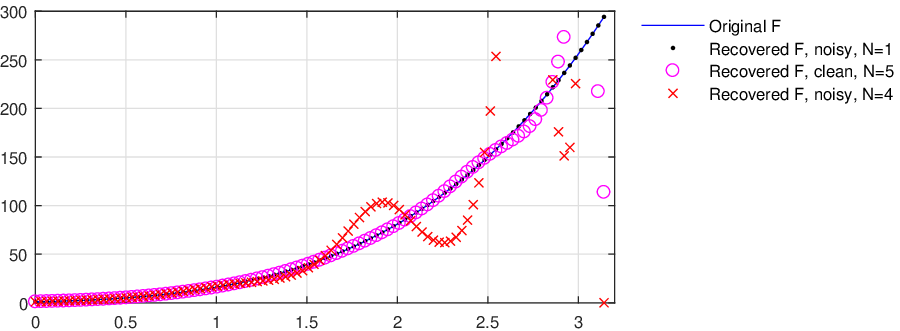}
\caption{
The original rod cross section area $F$ from Example  \ref{ExParabolic} together with the recovered ones using different values of $N_1=N_2=N$ in \eqref{equations 1} and clean or noisy data. On this figure we illustrate the importance of restricting the number of unknowns taken on Step \ref{AlgStep3}. The potential was better recovered from the noisy data by taking smaller values of $N_1=N_2=1$ than from the clean data with larger values  of $N_1$ and $N_2$.
}\label{Ex1FigPotential}
\end{figure}

\begin{figure}[htb!]
\centering
\includegraphics[bb=0 0 432 158, width=6in,height=2.2in]
{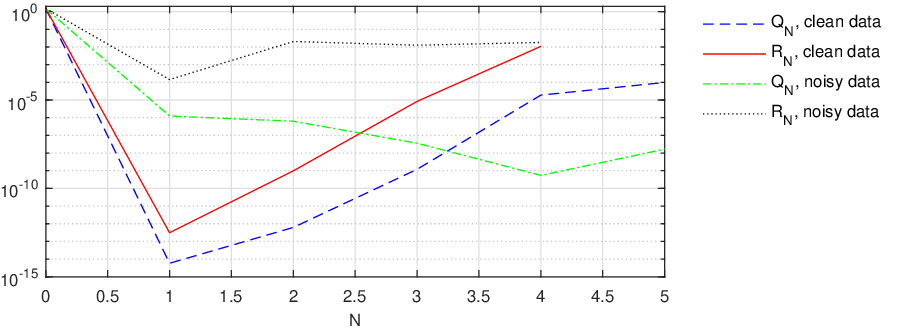}
\caption{Illustration to the work of Step \ref{AlgStep3} of the proposed algorithm for the inverse problem in Example \ref{ExParabolic}. Values of the functions $Q_N$ and $R_N$ given by \eqref{FunctionalOrig} and \eqref{Functional opt} are presented for both clean and noisy data}
\label{Ex1FigQNRN}
\end{figure}

\begin{figure}[htb!]
\centering
\includegraphics[bb=0 0 432 158, width=6in,height=2.2in]
{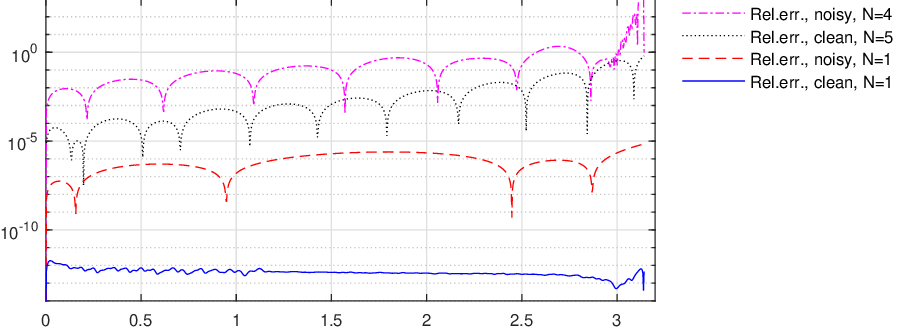}
\caption{Relative errors of the recovered rod cross section area $F$ from Example  \ref{ExParabolic} using different values of $N_1=N_2=N$ in \eqref{equations 1} and clean or noisy data.}
\label{Ex1FigRelErr}
\end{figure}

For the numerical experiment we took
\[
a=1,\quad b=1, \quad p=2,\quad E=3, \quad r=4
\]
and for the frequency, we took 12 equally spaced points, the first being $\omega_1=1$ and the last $\omega_{12}=2$. For the inverse problem we considered two sets of values. One being the exact values $\tilde f(\omega_k)$, $1\le k\le 12$ (referred to as ``clean data''), and the second one was obtained by adding noise using the formula
\[
\tilde f_\xi(\omega) = \tilde f(\omega) \cdot (1 + 10^{-6}\xi),
\]
here $\xi$ is the random value having standard normal distribution.

Even for the clean data, taking $N_1=N_2=5$ in \eqref{equations 1}, i.e., considering the number of unknowns to be equal to the number of equations, lead to a large error in the recovered shape $F$, as can be seen on Figures \ref{Ex1FigPotential} and \ref{Ex1FigRelErr}.

The choice of the parameters $N_1$ and $N_2$ based on minimizing the discrepancy \eqref{FunctionalOrig} results in $N_1=N_2=1$ for the clean data, but in $N_1=N_2=4$ for the noisy data, see Figure \ref{Ex1FigQNRN}. And for this latter choice, again, the coefficients $\{g_2^4,\ldots,g_4^4, s_2^4,\ldots, s_4^4\}$ obtained from \eqref{Eq argmin} are rather large, contrary to being equal to zero as was mentioned earlier. As a result, large error in the recovered shape $F$, see Figures \ref{Ex1FigPotential} and \ref{Ex1FigRelErr}.

The use of the function $R_N$ given by \eqref{Functional opt} gives us $N_1=N_2=1$ in both cases and the shape $F$ can be recovered with the relative error less than $2\cdot 10^{-12}$ for the clean data and less than $7\cdot 10^{-6}$ for the noisy data.
\end{example}

\begin{example}
\label{ExExp}
Consider a rod having exponential cross section area, i.e.,
\[
F(x) = \exp\left(2(a+bx)\right),
\]
where $a$ and $b$ are some constants. For this cross section profile we have
\[
a(x) = \exp(a+bx),\qquad q(x) = b^2,\qquad h = b.
\]
The solutions $\varphi$ and $S$ are known,
\begin{equation*}
    \varphi(\rho,x)=\cos \left(\sqrt{\rho^2 - b^2} x\right) + h \frac{\sin \left(\sqrt{\rho^2 - b^2} x\right)}{\sqrt{\rho^2 - b^2}},\qquad
    S(\rho,x)= \frac{\sin \left(\sqrt{\rho^2 - b^2} x\right)}{\sqrt{\rho^2 - b^2}},
\end{equation*}
hence we can compute $f(\rho)$ explicitly.

\begin{figure}[htb!]
\centering
\includegraphics[bb=0 0 432 158, width=6in,height=2.2in]
{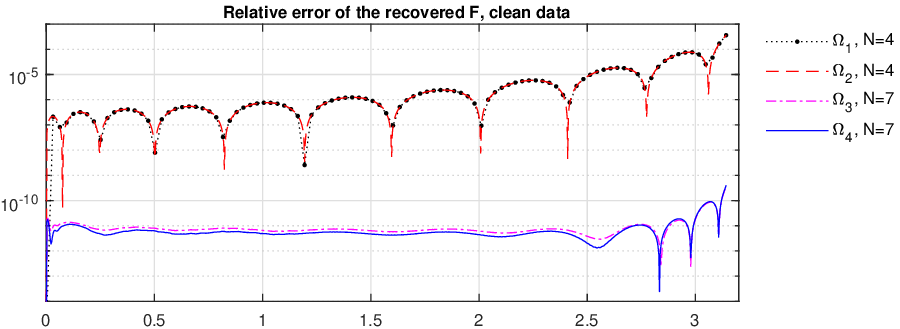}
\caption{
Relative errors of the recovered rod cross section area $F$ from Example  \ref{ExExp} using ``clean'' values $f(\omega)$ given in different sets of points $\Omega_j$, $j\in\{1,2,3,4\}$. The parameter $N$ for each plot indicates the value of $N_1=N_2=N$ taken in \eqref{equations 1} and obtained by minimizing the function $R_N$ \eqref{Functional opt}.
}\label{Ex2RelErrClean}
\end{figure}

\begin{figure}[htb!]
\centering
\includegraphics[bb=0 0 432 158, width=6in,height=2.2in]
{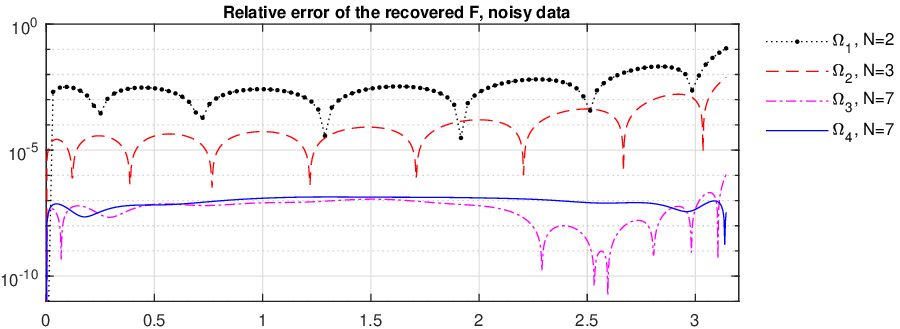}
\caption{Same as Figure \ref{Ex2RelErrClean}, but ``noisy'' data was used.}
\label{Ex2RelErrNoisy}
\end{figure}

\begin{figure}[htb!]
\centering
\includegraphics[bb=0 0 432 158, width=6in,height=2.2in]
{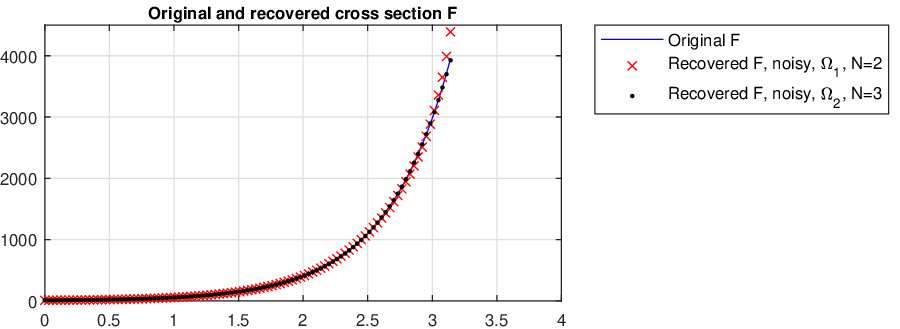}
\caption{The original cross section profile $F$ from Example \ref{ExExp} together with the recovered ones.}
\label{Ex2Profile}
\end{figure}

For the numerical experiment we took
\[
a=1,\quad b=1, \quad p=2,\quad E=3, \quad r=4
\]
and for the frequency, we considered $\Omega_1=\{1 + k/10, 0\le k\le 20\}$,
$\Omega_2=\{1 + k/40, 0\le k\le 80\}$, $\Omega_3=\{1 + 2k/5, 0\le k\le 20\}$ and $\Omega_4=\{1 + k/10, 0\le k\le 80\}$. I.e., $\Omega_1$ and $\Omega_2$ being 21 and 81 equally spaced points from the segment $[1,3]$, respectively, while $\Omega_3$ and $\Omega_4$ being 21 and 81 equally spaced points from the segment $[1,9]$.

As in Example \ref{ExParabolic}, we considered ``clean'' and ``noisy'' data, this time the strength of noise being $10^{-7}$.

On Figures \ref{Ex2RelErrClean} and \ref{Ex2RelErrNoisy} we present the relative error of the recovered rod cross section area. As one can see, the performance of the proposed method greatly depends on the interval where the original data is given, while the number of data points is not so important for the clean data. For the noisy data, having the data given at a larger number of points (even located in the same segment) improves precision and stability of the proposed method. On Figure \ref{Ex2Profile} one can see the original profile and the recovered ones. The visual difference can be observed for only one case, $\Omega_1$ set of frequencies and noisy data.
\end{example}

In the following examples we consider rods whose cross section areas are constant except two impurities. And the problem consists in determining the location (and shape) of these impurities. Similar problems were considered in \cite{WuFricke1990}.

\begin{example}\label{Ex2Holes}
For the first example we chose
\[
a(x) = \begin{cases}
1+\frac{1}{4}\left(1+\cos\left(6(x-\pi/3)\right)\right), & x\in(\pi/4,\pi/2),\\
1-\frac{1}{5}\left(1+\cos\left(8(x-3\pi/4)\right)\right), & x\in(5\pi/8,7\pi/8),\\
1, & \text{otherwise}.
\end{cases}
\]
The function $\tilde f(\omega)$ was computed on a uniform mesh of 101 points covering the segment $[0.1,20]$. In the variable $\rho$ this segment converts to $[0.0866, 17.3205]$. This problem does not possess exact solution so the direct data $f(\rho)$ was computed numerically using the method from \cite{KNT}.

In this example more coefficients $g_n$ and $s_n$ are required for good approximation of the solutions $\varphi$ and $S$. The minimum of the functional $Q_N$ given by \eqref{FunctionalOrig} was achieved for $N=19$, see Figure \ref{Ex3Figs}. However, the coefficients $g_0^{19},\ldots, g_{19}^{19}$ computed from \eqref{Eq argmin} were large by absolute value and quite distant from the coefficients obtained if one takes $N=18$. This indicates possible large error of the approximate solution $\varphi_{19}(\rho,\pi)$, especially for the values $\rho$ lying outside the given interval, i.e., for $\rho> 17.3$. We verified this by checking errors of the computed eigenvalues $\mu_n$, see Figure \ref{Ex3Figs}. We additionally considered smaller $N=9$ (where the discrepancy in the overdetermined system \eqref{equations 1}--\eqref{equations 3} is larger, but numerical errors are expected to be less significant), and $N=42$ (obtained as the minimum of the functional $R_N$). As one can see, the recovered coefficients are close to the original ones only for $N=42$. Nevertheless, the eigenvalues are recovered quite well inside the given interval of $\rho$ (i.e., for $\rho\le 17.3$), however the errors increase outside of this interval, especially for $N=19$, leading to poorly recovered potential. For this reason we recommend to use the functional $R_N$.

\begin{figure}[htb!]
\centering
\begin{tabular}{cc}
\includegraphics[bb=0 0 216 173, width=3in,height=2.4in]
{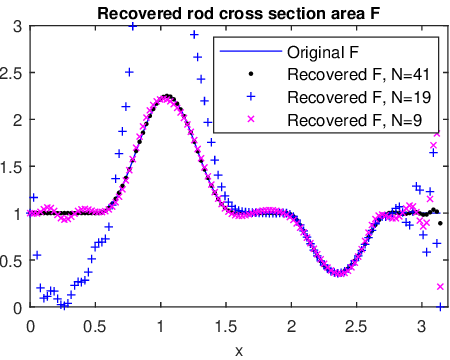} &
\includegraphics[bb=0 0 216 173, width=3in,height=2.4in]
{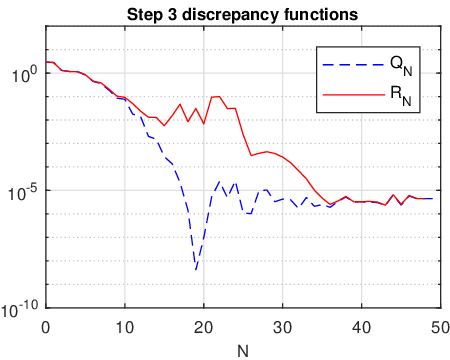}
\\
\includegraphics[bb=0 0 216 173, width=3in,height=2.4in]
{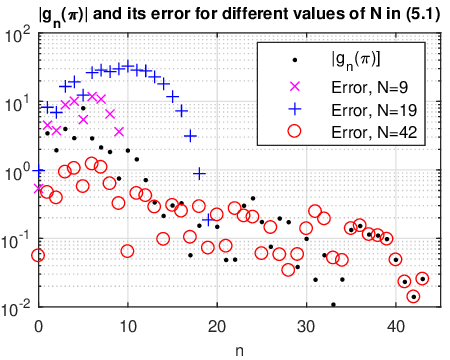} &
\includegraphics[bb=0 0 216 173, width=3in,height=2.4in]
{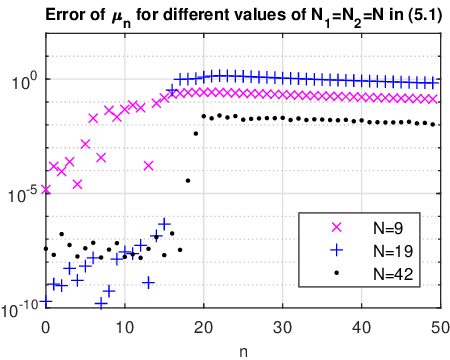}
\end{tabular}
\caption{Illustration to the work of the proposed algorithm for the inverse problem in Example \ref{Ex2Holes}. In all the plots parameter $N$ corresponds to the values $N_1=N_2=N$ taken on Step \ref{AlgStep3}.
Top left plot: original and recovered rod cross section area $F$.
Top right plot: values of the functions $Q_N$ and $R_N$ given by \eqref{FunctionalOrig} and \eqref{Functional opt}. Bottom left plot: absolute errors in the recovered coefficients $g_n(\pi)$ on Step \ref{AlgStep3}. Bottom right plot: accuracy of the eigenvalues $\mu_n$ found on the Step \ref{AlgStep4}.
\label{Ex3Figs}}
\end{figure}

\end{example}

\begin{example}\label{Ex2HolesSm}
For the second example we chose
\[
a(x) = \begin{cases}
1+\frac{1}{10}\exp\left(1- \frac{\pi^2}{(5\pi-12x)(12x-3\pi)}\right), & x\in(\pi/4,5\pi/12),\\
1-\frac{1}{15}\exp\left(1- \frac{\pi^2}{(31\pi-40x)(40x-29\pi)}\right), & x\in(29\pi/40,31\pi/40),\\
1, & \text{otherwise}.
\end{cases}
\]
The function $\tilde f(\omega)$ was computed on a uniform mesh of 201 points covering the segment $[0.1,50]$. In the variable $\rho$ this segment converts to $[0.0866, 43.30]$. This problem does not possess exact solution so the direct data $f(\rho)$ was computed numerically using the method from \cite{KNT}. We also considered two additional datasets taking the first 41 and 81 points respectively.

On Figure \ref{Ex4Fig} we show the recovered rod cross section area $F$ for these three datasets. One can appreciate that the accuracy increases with additional data available. The functional $R_N$ was utilized to determine optimal value of $N_1=N_2=N$ for the Step \ref{AlgStep3}. $N=18$ was used for 41-point dataset, $N=38$ for  81-point dataset and $N=85$ for 201-point dataset.

\begin{figure}[htb!]
\centering
\includegraphics[bb=0 0 360 216, width=5in,height=3in]
{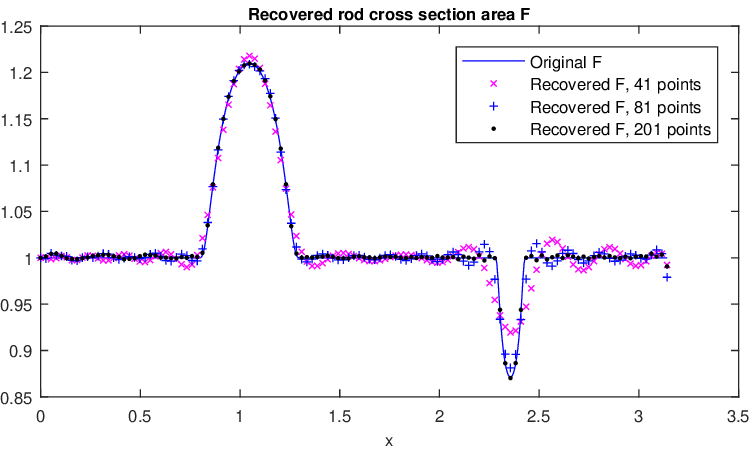}
\caption{Illustration to the work of the proposed algorithm for the inverse problem in Example \ref{Ex2HolesSm}. Different number of data points used: 41, 81 and 201 points $\omega_l$.}
\label{Ex4Fig}
\end{figure}
\end{example}

\section*{Funding information}
CONAHCYT, Mexico, grant ``Ciencia de Frontera'' FORDECYT--PRONACES/61517/2020.

\section*{Data availability}
The data that support the findings of this study are available upon reasonable request.

\section*{Conflict of interest}
This work does not have any conflict of interest.

%

\end{document}